\documentclass[
framed_theorems,
framed_proofs,
print,
]{OCPreprint}

\usepackage{enumitem}

\renewcommand\d{\mathop{}\!\mathrm{d}}

\title[Slater conditions without interior points]
{Slater conditions without interior points
for programs in Lebesgue spaces with pointwise bounds and finitely many constraints}
\author{Gerd Wachsmuth\footnote{%
		Brandenburgische Technische Universität Cottbus--Senftenberg,
		Institute of Mathematics,
		03046 Cottbus,
		Germany,
		\email{wachsmuth@b-tu.de},
		\url{https://www.b-tu.de/fg-optimale-steuerung/team/prof-gerd-wachsmuth},
		ORCID: 0000-0002-3098-1503%
	}
}
{
	\makeatletter
	\def\and{ and }
	\def\footnote#1{}
	\hypersetup{
		pdftitle={\@title},
		pdfauthor={\@author}
	}
	\makeatother
}

\begin{document}
%%fakesection: Title
\maketitle
\begin{abstract}
	We consider optimization problems in Lebesgue spaces
	with pointwise box constraints and finitely many additional
	linear constraints.
	We prove that the existence of a Slater point
	which lies strictly between the pointwise bounds
	and which satisfies the linear constraints
	is sufficient for the existence of Lagrange multipliers.
	Surprisingly, the Slater point is also necessary
	for the existence of Lagrange multipliers in a certain sense.
	We also demonstrate how to handle additional
	finitely many nonlinear constraints.
\end{abstract}


\begin{msc}
	\mscLink{49K27},
	\mscLink{90C46}
\end{msc}

\section{Introduction}
\label{sec:intro}
We consider the problem
\begin{equation*}
	\begin{aligned}
		\text{Minimize} \quad & f(x)\\
		\text{s.t.} \quad
		&
		\begin{aligned}[t]
			x_a \le x &\le x_b && \text{$\mu$-a.e.\ on $\Omega$}, &
			\dual{g_i}{x}_{L^p(\mu)} &\le a_i && \forall i = 1,\ldots, n,\\
			&&&&
			\dual{h_j}{x}_{L^p(\mu)} &= b_j && \forall j = 1,\ldots, m.
		\end{aligned}
	\end{aligned}
\end{equation*}
In this problem, we have pointwise bounds $x_a$, $x_b$,
and finitely many linear constraints.
We are interested in the existence of Lagrange multipliers.
This is a delicate question,
since (under typical assumptions on the measure $\mu$)
the set
\begin{equation*}
	K := \set{x \in L^p(\mu) \given x_a \le x \le x_b \; \text{$\mu$-a.e.}}
\end{equation*}
has empty interior
unless $p = \infty$.
Moreover,
if $p = \infty$,
this would lead initially to multipliers
in the space $L^\infty(\mu)\dualspace$,
but often their regularity can be improved.
We mention that the above problem
is also interesting from an application point of view,
see, e.g.,
\cite{BonnansZidani1999,CasasTroeltzsch2002,BuchheimGrueteringMeyer2022}.

We say that $\hat x \in L^p(\mu)$ is a Slater point,
if the conditions
\begin{align*}
	x_a < \hat x &< x_b \quad \text{$\mu$-a.e.\ on $\Omega$},
	&
	\dual{g_i}{\hat x}_{L^p(\mu)} &\le a_i \quad \forall i = 1,\ldots, n,\\
	&&
	\dual{h_j}{\hat x}_{L^p(\mu)} &= b_j \quad \forall j = 1,\ldots, m
\end{align*}
are satisfied.
Note that we only require that $\hat x$ lies strictly between the pointwise bounds
and this does not imply that $\hat x$ is an interior point of $K$.
Our main result \cref{thm:foc_linear}
shows the existence of Lagrange multipliers for locally optimal solutions $\bar x$
of the above problem,
provided that a Slater point exists.
Under certain circumstances,
we are even able to prove that a Slater point is necessary for the existence of Lagrange
multipliers (see \cref{sec:slater_nec}),
which is rather surprising.
We also address the case of nonlinear inequality constraints
in \cref{sec:nonlinear}.

We briefly describe our idea to prove the existence of multipliers.
We denote by $P$ the set of points satisfying the finite-dimensional constraints,
i.e., the feasible set of our problem is given by $K \cap P$.
For the existence of Lagrange multipliers,
we have to characterize the normal cone of $K \cap P$.
By using the notion of $n$-polyhedricity from \cite{Wachsmuth2016:2},
we automatically get the formula
$\TT_{K \cap P}(x) = \TT_K(x) \cap \TT_P(x)$
for the tangent cone.
By dualizing this formula,
we get
$\NN_{K \cap P}(x) = \NN_K(x) + \NN_P(x)$
provided that the set on the right-hand side is closed
(or weak-$\star$ closed, see \cref{lem:polar_of_intersection} below).
This closedness can be shown by utilizing the Slater point $\hat x$.

We mention that a similar result is obtained in the rather unknown work
\cite{Troeltzsch1977}.
Therein, the existence of Lagrange multipliers is proved
using a similar notion of Slater points.
In contrast to our approach,
the main tool is duality for linear programs.
Under rather weak conditions, it is shown that no duality gap occurs
and the Slater point is used to prove the existence of a dual solution.
Moreover,
the problem is addressed in more abstract spaces.
However,
our approach has certain advantages.
We are able to discuss also $p \in \set{1,\infty}$,
whereas
\cite{Troeltzsch1977}
requires reflexive spaces.
We further show necessity of Slater points for the existence of multipliers.
We address inequality constraints
whereas the approach of
\cite{Troeltzsch1977}
is only given for equality constraints.

\section{Preliminaries}
\label{sec:prelim}

Throughout the paper,
we use the following setting.
\begin{assumption}[Standing assumption]
	\label{asm:standing}
	The triple
	$(\Omega, \Sigma, \mu)$
	is a $\sigma$-finite
	measure space
	and
	$p \in [1,\infty]$ is fixed.
	By
	$p' \in [1,\infty]$
	we denote its conjugate exponent,
	i.e., $1/p + 1/p' = 1$ with $1/\infty = 0$.
\end{assumption}

We fix some notation.
The set of nonnegative integers is denoted by $\N$.
The set of positive, real numbers is $\R^+$.
By $L^q(\mu)$ we denote the usual Lebesgue spaces.
Recall that $L^q(\mu)\dualspace = L^{q'}(\mu)$
for $q \in [1,\infty)$,
see \cite[243G and 244K]{Fremlin2003}.

Finally, we mention that we use $\set{x_a = x}$
with functions $x_a, x \colon \Omega \to \R$
as an abbreviation for the set
$\set{\omega \in \Omega \given x_a(\omega) = x(\omega)}$.
Similar abbreviations
are also used for different relation symbols.

\subsection{Tangent cones and normal cones}
\label{subsec:cones}
We recall some preliminary results
concerning tangent cones and normal cones.
Further information can be found in \cite{BonnansShapiro2000}.
\begin{definition}
	\label{def:tangent_cone}
	Let $C \subset X$ be a closed subset of the Banach space $X$.
	For $x \in C$, we define the radial cone
	\begin{equation*}
		\RR_C(x)
		:=
		\set{
			d \in X
			\given
			\exists t_0 > 0 :
			\forall t \in (0,t_0) :
			x + t d \in C
		}
	\end{equation*}
	and the (Bouligand) tangent cone
	\begin{equation*}
		\TT_C(x)
		:=
		\set{
			d \in X
			\given
			\exists\seq{x_k}\subset C, \seq{t_k}\subset\R^+ :
			x_k\to x
			,\;
			t_k\searrow 0
			,\;
			(x_k - x)/t_k\to d
		}
		.
	\end{equation*}
	For $x \in X \setminus C$, we set $\RR_C(x) := \TT_C(x) := \emptyset$.
\end{definition}
In the case that $C$ is convex, the tangent cone $\TT_C(x)$
coincides with the closure of the radial cone $\RR_C(x)$.

For box constraints in $L^p(\mu)$, $p \in [1,\infty)$,
the following formula for the tangent cone is well known.
We also provide estimates for the tangent cone in $L^\infty(\mu)$.
\begin{lemma}
	\label{lem:tangent_cone}
	Let measurable functions $x_a, x_b \colon \Omega \to [-\infty,\infty]$ be given
	such that
	\begin{equation*}
		K
		:=
		\set{ x \in L^p(\mu) \given x_a \le x \le x_b \; \text{$\mu$-a.e.\ on $\Omega$}}
	\end{equation*}
	is nonempty
	and we fix $x \in K$.
	The estimates
	\begin{subequations}
		\label{eq:tangent}
		\begin{align}
			\label{eq:tangent:1}
			\TT_K(x)
			&\supset
			\set*{ h \in L^\infty(\mu) \cap L^p(\mu) \given
				\begin{aligned}
					h &\ge 0 \; \text{$\mu$-a.e.\ on $\set{x_a + \varepsilon \ge x}$}, \\
					h &\le 0 \; \text{$\mu$-a.e.\ on $\set{x_b - \varepsilon \le x}$}
				\end{aligned}
			},
			\\
			\label{eq:tangent:2}
			\TT_K(x)
			&\subset
			\set{ h \in L^p(\mu) \given
				h \ge 0 \; \text{$\mu$-a.e.\ on $\set{x_a = x}$},
				\;
				h \le 0 \; \text{$\mu$-a.e.\ on $\set{x_b = x}$}
			}
		\end{align}
	\end{subequations}
	hold for all $\varepsilon > 0$.
	For $p \in [1,\infty)$, we even have
	\begin{equation}
		\label{eq:tangent}
		\TT_K(x)
		=
		\set{ h \in L^p(\mu) \given
			h \ge 0 \; \text{$\mu$-a.e.\ on $\set{x_a = x}$},
			\;
			h \le 0 \; \text{$\mu$-a.e.\ on $\set{x_b = x}$}
		}.
	\end{equation}
\end{lemma}
\begin{proof}
	The estimates are straightforward to verify.

	Now, let $p < \infty$ and let $h$ belong to the right-hand side of \eqref{eq:tangent}.
	We define the sequence $h_k \in L^p(\mu)$ via
	\begin{equation*}
		h_k(\omega)
		:=
		\begin{cases}
			0 & \text{if } x(\omega) \in (x_a(\omega), x_a(\omega) + 1/k) \cup ( x_b(\omega) - 1/k, x_b(\omega)) \text{ or } \abs{h(\omega)} > k, \\
			h(\omega) & \text{else}.
		\end{cases}
	\end{equation*}
	Then, it is easy to check that $x_k = x + h_k/k^2 \in K$.
	Moreover, we have $h_k \to h$ pointwise and, due to $h \in L^p(\mu)$,
	the dominated convergence theorem gives $h_k \to h$ in $L^p(\mu)$.
	This yields $h \in \TT_K(x)$.
\end{proof}
In case $p = \infty$, we do not have $h_k \to h$ in $L^\infty(\mu)$,
but one can show $h_k \weaklystar h$ in $L^\infty(\mu)$.
Thus, an analogue formula to \eqref{eq:tangent}
holds, if we define the tangent cone in $L^\infty(\mu)$
by using the weak-$\star$ topology.

Our next goal is the definition of the normal cone
which utilizes the concept of the polar cone.
When we are working in the space $L^\infty(\mu)$,
we do not want to arrive in the cumbersome space
$L^\infty(\mu)\dualspace$,
but in the nice predual space $L^1(\mu)$.
This is already reflected in the definition of the polar cone.
\begin{definition}
	\label{def:polar}
	Let $C \subset X$ be a subset of a Banach space $X$.
	We define
	\begin{equation*}
		C\polar
		:=
		\set{
			y \in L^1(\mu)
			\given
			\forall x \in C :
			\dual{y}{x}_{L^\infty(\mu)} \le 0
		}
	\end{equation*}
	in case $X = L^\infty(\mu)$
	and
	\begin{equation*}
		C\polar
		:=
		\set{
			y \in X\dualspace
			\given
			\forall x \in C :
			\dual{y}{x}_{X} \le 0
		}
	\end{equation*}
	otherwise.
\end{definition}
Another perspective
is to equip the pair of spaces $L^p(\mu)$ and $L^{p'}(\mu)$
with compatible topologies.
That is, we can use the norm topologies
in $L^p(\mu)$, $p \in [1,\infty)$,
but
$L^\infty(\mu)$
has to be equipped with the weak-$\star$ topology.
In this sense,
\cref{def:polar} is compatible with, e.g.,
\cite[Section~2.1.4]{BonnansShapiro2000}.

As usual,
the normal cone is defined as the polar of the tangent cone.
\begin{definition}
	\label{def:normal}
	Let $C \subset X$ be a closed subset of a Banach space $X$.
	For $x \in C$ we define
	\begin{equation*}
		\NN_C(x) := \TT_C(x)\polar
	\end{equation*}
	and
	$\NN_C(x) := \emptyset$ for $x \in X \setminus C$.
\end{definition}
Note that the polar of $\TT_C(x)$ is defined as in \cref{def:polar}.
In particular,
this is the usual definition of the normal cone,
unless $X = L^\infty(\mu)$.
In this case case,
we define the normal cone to stay in
the predual space $L^1(\mu)$ of $L^\infty(\mu)$.

\begin{lemma}
	\label{lem:normal_cone}
	Let $K$ and $x$ be given as in \cref{lem:tangent_cone}.
	Then,
	\begin{equation*}
		\NN_K(x)
		=
		\set*{
			\zeta \in L^{p'}(\mu)
			\given
			\begin{aligned}
				& \zeta \le 0 \; \text{$\mu$-a.e.\ on $\set{x_a = x < x_b}$},
				\\
				& \zeta = 0 \; \text{$\mu$-a.e.\ on $\set{x_a < x < x_b}$},
				\\
				& \zeta \ge 0 \; \text{$\mu$-a.e.\ on $\set{x_a < x = x_b}$}
			\end{aligned}
		}
		.
	\end{equation*}
\end{lemma}
\begin{proof}
	The inclusion ``$\supset$'' follows easily from the estimate \eqref{eq:tangent:2}.

	Now, let $\zeta \in \NN_K(x)$ be given.
	Let $\varepsilon > 0$ be given and let $M$ be an arbitrary subset of
	$\set{ x_a = x \le x_b - \varepsilon} \cap \set{\zeta > 0}$
	with $\mu(M) < \infty$.
	By \eqref{eq:tangent:1}, its characteristic function $\xi_M$
	belongs to $\TT_K(x)$, thus
	\begin{equation*}
		0 \ge \dual{\chi_M}{\zeta} = \int_M \zeta \d\mu = \int_M \abs{\zeta} \d\mu.
	\end{equation*}
	Due to $\zeta > 0$ $\mu$-a.e.\ on $M$, this implies $\mu(M) = 0$.
	Consequently,
	$\set{ x_a = x \le x_b - \varepsilon} \cap \set{\zeta > 0}$
	is a $\mu$-null set.
	By considering $\varepsilon \to 0$, also
	$\set{ x_a = x < x_b} \cap \set{\zeta > 0}$
	is a $\mu$-null set
	and, hence,
	$\zeta \le 0$ $\mu$-a.e.\ on $\set{ x_a = x < x_b}$.
	The other sign conditions on $\zeta$ can be shown analogously.
\end{proof}

For later reference,
we provide the following result
for the polar of an intersection of tangent cones.
\begin{lemma}
	\label{lem:polar_of_intersection}
	Let $C, P \subset L^p(\mu)$ be closed and convex
	and fix $x \in C \cap P$.
	Then,
	\begin{equation*}
		\parens*{
			\TT_C(x) \cap \TT_P(x)
		}\polar
		=
		\cl_\star\parens*{
			\NN_C(x) + \NN_P(x)
		}
		.
	\end{equation*}
	Here, $\cl_\star(A)$ for $A \subset L^{p'}(\mu)$
	denotes
	the closure of $A$ w.r.t.\ the norm convergence in case $p \in (1,\infty]$
	and
	the closure of $A$ w.r.t.\ the weak-$\star$ topology of $L^\infty(\mu)$ in case $p = 1$.
\end{lemma}
\begin{proof}
	This follows from
	\cite[(2.32)]{BonnansShapiro2000}.
\end{proof}

Finally,
we provide simple formulas for sets defined by linear restrictions.
\begin{lemma}
	\label{lem:finite_restrictions}
	Let $n,m \in \N$,
	$g_i, h_j \in L^{p'}(\mu)$ and
	$a_i, b_j \in \R$, $1 \le i \le n$, $1 \le j \le m$
	be given.
	We define
	\begin{equation*}
		P
		:=
		\set*{ x \in L^p(\mu) \given
			\begin{aligned}
				\dual{g_i}{x}_{L^p(\mu)} &\le a_i && \forall i = 1,\ldots, n,\\
				\dual{h_j}{x}_{L^p(\mu)} &= b_j && \forall j = 1,\ldots, m
			\end{aligned}
		}
		.
	\end{equation*}
	For all $x \in P$ we have
	\begin{align*}
		\TT_P(x)
		&=
		\set*{ y \in L^p(\mu) \given
			\begin{aligned}
				\dual{g_i}{y}_{L^p(\mu)} &\le 0 && \forall i \in \AA(x),\\
				\dual{h_j}{y}_{L^p(\mu)} &= 0 && \forall j = 1,\ldots, m
			\end{aligned}
		},
		\\
		\NN_P(x)
		&=
		\set[\Bigg]{
			\sum_{i \in \AA(x)} \alpha_i g_i
			+
			\sum_{j = 1}^m \beta_j h_j
			\in
			L^{p'}(\mu)
			\given
			\alpha_i \ge 0,\, i \in \AA(x),\;
			\beta_j \in \R,\, j = 1,\ldots,m
		},
	\end{align*}
	where
	\begin{equation*}
		\AA(x) := \set[\big]{ i \in \set{1,\ldots,n} \given \dual{g_i}{x}_{L^p(\mu)} = a_i }
	\end{equation*}
	is the active set.
\end{lemma}
\begin{proof}
	The formula for the tangent cone follows from a simple calculation.
	The normal cone is given in \cite[Proposition~2.42]{BonnansShapiro2000}.
\end{proof}

\subsection{Polyhedric sets}
\label{subsec:poly}
In this section, we briefly review some of the theory of polyhedric sets.
More details can be found in \cite{Wachsmuth2016:2}.

Polyhedric sets are a suitable generalization of polytopes
to infinite dimensional spaces.
In many regards,
they seem to behave similarly to polytopes,
in particular if they are in interaction with
sets defined via finitely many continuous and linear inequalities.
They were introduced
in the seminal works \cite{Mignot1976,Haraux1977}.
In \cite{Wachsmuth2016:2}
it was shown that ``almost all'' of the polyhedric sets
are even $n$-polyhedric for all $n \in \N$.

\begin{definition}
	\label{def:higher_order_polyhedricity}
	Let $C \subset X$ be a closed, convex subset of the Banach space $X$
	and let $n \in \N$ be given.
	We call $C$ $n$-polyhedric at $x \in C$,
	if
	\begin{equation}
		\label{eq:npol}
		\TT_C(x) \cap \bigcap_{i = 1}^n \xi_i\anni
		=
		\cl\braces[\Big]{
			\RR_C(x) \cap \bigcap_{i = 1}^n \xi_i\anni
		}
		\qquad
		\forall \xi_1,\ldots,\xi_n \in X\dualspace
	\end{equation}
	holds,
	where $\xi\anni := \set{y \in X \given \dual{\xi}{y}_X = 0}$.
	Moreover, $C$ is called $n$-polyhedric,
	if it is $n$-polyhedric at all $x \in C$.
\end{definition}
The classical case of polyhedricity corresponds to
the case $n = 1$ in \cref{def:higher_order_polyhedricity}.

A main result of \cite{Wachsmuth2016:2} is
that sets defined by lower and upper bounds
are $n$-polyhedric for all $n \in \N$
if the underlying Banach space
possesses a certain lattice structure.
This, in particular, applies
to our situation
and yields the following result, cf.\ \cite[Example~4.21(1)]{Wachsmuth2016:2}.
\begin{theorem}
	\label{thm:K_n_poly}
	The set $K$ from \cref{lem:tangent_cone}
	is $n$-polyhedric for all $n \in \N$.
\end{theorem}

This $n$-polyhedricity
yields an easy formula for the tangent cone and the normal cone
of the intersection with finitely many linear restrictions.
\begin{lemma}
	\label{lem:poly_cones}
	Let subsets $K,P \subset L^p(\mu)$ be given.
	We assume that $K$ is $n$-polyhedric for all $n \in \N$
	and that $P$ is defined by finitely many continuous and linear
	equalities and inequalities.
	Then,
	\begin{align*}
		\TT_{K \cap P}(x)
		&=
		\TT_K(x) \cap \TT_P(x)
		,
		&
		\NN_{K \cap P}(x)
		&=
		\cl_\star\parens*{
			\NN_K(x) + \NN_P(x)
		}
	\end{align*}
	for all $x \in K \cap P$.
	Here, $\cl_\star(A)$ for $A \subset L^{p'}(\mu)$
	denotes
	the closure of $A$ w.r.t.\ the norm convergence in case $p \in (1,\infty]$
	and
	the closure of $A$ w.r.t.\ the weak-$\star$ topology of $L^\infty(\mu)$ in case $p = 1$.
\end{lemma}
\begin{proof}
	The formula for the tangent cone follows
	from \cite[Lemma~4.4]{Wachsmuth2016:2}.
	Consequently, \cref{lem:polar_of_intersection}
	yields the identity for the normal cone.
\end{proof}
Recall that $\TT_P(x)$ and $\NN_P(x)$ are given in \cref{lem:finite_restrictions}.

\section{Linear Constraints}
\label{sec:linear}
In this section, we focus on a problem
with linear constraints
of the form
\begin{equation}
	\label{eq:linear}
	\begin{aligned}
		\text{Minimize} \quad & f(x)\\
		\text{w.r.t.} \quad & x \in L^p(\mu) \\
		\text{s.t.} \quad
		&
		\begin{aligned}[t]
			x_a \le x &\le x_b && \text{$\mu$-a.e.\ on $\Omega$}, \\
			\dual{g_i}{x}_{L^p(\mu)} &\le a_i && \forall i = 1,\ldots, n,\\
			\dual{h_j}{x}_{L^p(\mu)} &= b_j && \forall j = 1,\ldots, m.
		\end{aligned}
	\end{aligned}
\end{equation}
Here,
$f \colon L^p(\mu) \to \R$
is assumed to be Fréchet differentiable,
$x_a, x_b \colon \Omega \to [-\infty,\infty]$
are measurable,
$n,m \in \N$,
$g_i, h_j \in L^{p'}(\mu)$
and
$a_i, b_j \in \R$.
Moreover, in case $p = \infty$,
we require
$f'(x) \in L^{p'}(\mu) = L^1(\mu)$
for all feasible points $x$ of \eqref{eq:linear}.

We denote by $F \subset L^p(\mu)$ the feasible set of \eqref{eq:linear}.
By standard arguments, every local minimizer $\bar x$ satisfies
$\dual{f'(\bar x)}{d}_{L^p(\mu)} \ge 0$
for all $d \in \TT_F(\bar x)$
and, consequently,
\begin{equation}
	\label{eq:opt_con_lin}
	\mathopen % Otherwise ugly spacing.
	-f'(\bar x)
	\in
	\NN_F(\bar x).
\end{equation}
It remains to
evaluate the normal cone.
To this end,
we write the feasible set
as $F = K \cap P$,
where
\begin{align*}
	K
	&:=
	\set{ x \in L^p(\mu) \given x_a \le x \le x_b \; \text{$\mu$-a.e.\ on $\Omega$}}
	\\
	P
	&:=
	\set*{ x \in L^p(\mu) \given
		\begin{aligned}
			\dual{g_i}{x}_{L^p(\mu)} &\le a_i && \forall i = 1,\ldots, n,\\
			\dual{h_j}{x}_{L^p(\mu)} &= b_j && \forall j = 1,\ldots, m
		\end{aligned}
	}
	.
\end{align*}
Using the results on $n$-polyhedric sets recalled in \cref{subsec:poly},
we get
\begin{equation*}
	\NN_F(\bar x)
	=
	\NN_{K \cap P}( \bar x )
	=
	\cl_\star \parens*{
		\NN_K(\bar x)
		+
		\NN_P(\bar x)
	}.
\end{equation*}
However, the set
$\NN_K(\bar x) + \NN_P(\bar x)$
is, in general,
not (weak-$\star$)
closed,
see \cref{sec:counterex}.
Thus, we have to utilize a Slater point.

\begin{definition}
	\label{def:slater_point}
	A point $\hat x \in L^p(\mu)$
	is called Slater point of \eqref{eq:linear},
	if
	the conditions
	\begin{align*}
		x_a < \hat x &< x_b \quad \text{$\mu$-a.e.\ on $\Omega$},
		&
		\dual{g_i}{\hat x}_{L^p(\mu)} &\le a_i \quad \forall i = 1,\ldots, n,\\
		&&
		\dual{h_j}{\hat x}_{L^p(\mu)} &= b_j \quad \forall j = 1,\ldots, m
	\end{align*}
	are satisfied.
	We say that \eqref{eq:linear} satisfies
	the Slater condition
	if there exists a Slater point of \eqref{eq:linear}.
\end{definition}
In the next results,
we use the formulas
for $\NN_K(x)$ and $\NN_P(x)$
provided in \cref{lem:normal_cone,lem:finite_restrictions}.
\begin{lemma}
	\label{lem:slater}
	Suppose that $\hat x$ is a Slater point of \eqref{eq:linear}.
	Then, for every feasible point $x$ of \eqref{eq:linear},
	the set
	$\NN_K(x) + \NN_P(x) \subset L^{p'}(\mu)$
	is closed ($p \in (1,\infty]$)
	or weak-$\star$ closed ($p = 1$), respectively.
\end{lemma}
\begin{proof}
	We begin with the easier case $p > 1$.
	Suppose that
	the sequence $\seq{\xi_k} \subset \NN_K(x) + \NN_P(x)$
	is convergent in $L^{p'}(\mu)$,
	i.e.,
	\begin{equation*}
		\xi_k
		=
		\zeta_k + \sum_{i \in \AA(x)} \alpha_{k,i} g_i + \sum_{j = 1}^m \beta_{k,j} h_j
		\to
		\xi
	\end{equation*}
	in $L^{p'}(\mu)$,
	where $\zeta_k \in \NN_K(x)$, $\alpha_{k,i} \ge 0$ and $\beta_{k,j} \in \R$.
	This implies
	\begin{equation}
		\label{eq:slater_conv834}
		\dual{\xi_k}{\hat x - x}_{L^p(\mu)}
		=
		\dual{\zeta_k}{\hat x - x}_{L^p(\mu)}
		+
		\sum_{i \in \AA(x)} \alpha_{k,i} \dual{g_i}{\hat x - x}_{L^p(\mu)}
		+
		\sum_{j = 1}^m \beta_{k,j} \dual{h_j}{\hat x - x}_{L^p(\mu)}
		.
	\end{equation}
	Due to the feasibility of the Slater point,
	we have
	\begin{equation*}
	\begin{aligned}
		\dual{\zeta_k}{\hat x - x}_{L^p(\mu)} &\le 0, \\
		\dual{g_i}{\hat x - x}_{L^p(\mu)}
		&=
		\dual{g_i}{\hat x}_{L^p(\mu)} - a_i \le 0 &&\forall i \in \AA(x), \\
		\dual{h_j}{\hat x - x}_{L^p(\mu)} &= 0 &&\forall j = 1,\ldots, m.
	\end{aligned}
	\end{equation*}
	Thus, all addends on the right-hand side of \eqref{eq:slater_conv834} are nonpositive
	and this yields
	\begin{equation}
		\label{eq:slater_conv834_bounds}
		0
		\ge
		\dual{\zeta_k}{\hat x - x}_{L^p(\mu)}
		\ge
		\dual{\xi_k}{\hat x - x}_{L^p(\mu)}
		.
	\end{equation}
	Since the right-hand side is a convergent sequence in $\R$,
	we get the boundedness of the real-valued sequence
	$\seq{\dual{\zeta_k}{\hat x - x}_{L^p(\mu)}}_k$.
	Next, we choose a function $y \in L^p(\mu)$
	which is positive $\mu$-a.e.\ and satisfies
	$y = \abs{\hat x - x}$ a.e.\ on $\set{\hat x \ne x}$.
	Further, we define the measure $\hat\mu$ via $\hat \mu = y \mu$.
	Due to $y \in L^p(\mu)$, we have the continuous embedding
	$L^{p'}(\mu) \embeds L^1(\hat\mu)$.
	Since $\zeta_k$ vanishes $\mu$-a.e.\ on $\set{\hat x = x}$
	and $\zeta_k (\hat x - x) \le 0$ $\mu$-a.e.\ on $\Omega$,
	we get
	\begin{equation*}
		\norm{\zeta_k}_{L^1(\hat\mu)}
		=
		\int_\Omega \abs{\zeta_k} \d \hat\mu
		=
		\int_\Omega \abs{\zeta_k} y \d \mu
		=
		\int_\Omega \abs{\zeta_k} \abs{\hat x - x} \d \mu
		% =
		% -\int_\Omega \zeta_k (\hat x - x) \d\mu
		=
		-\dual{\zeta_k}{\hat x - x}_{L^p(\mu)}.
	\end{equation*}
	This shows that the sequence $\zeta_k$ is bounded in the space $L^1(\hat\mu)$.
	Since $\xi_k$ is bounded in $L^{p'}(\mu)$,
	it is also bounded in $L^1(\hat\mu)$.
	Thus, the sequence
	\begin{equation*}
		\seq{ \xi_k - \zeta_k }_k
		=
		\seq*{\sum_{i \in \AA(x)} \alpha_{k,i} g_i + \sum_{j = 1}^m \beta_{k,j} h_j}_k
	\end{equation*}
	is bounded in $L^1(\hat\mu)$.
	Now, we can invoke \cite[Proposition~2.41]{BonnansShapiro2000}
	and this implies the existence of bounded sequences
	$\seq{\hat\alpha_{k,i}}_k \subset [0,\infty)$,
	$\seq{\hat\beta_{k,j}}_k \subset \R$,
	such that
	\begin{equation*}
		\sum_{i \in \AA(x)} \alpha_{k,i} g_i + \sum_{j = 1}^m \beta_{k,j} h_j
		=
		\sum_{i \in \AA(x)} \hat\alpha_{k,i} g_i + \sum_{j = 1}^m \hat\beta_{k,j} h_j
	\end{equation*}
	for all $k$.
	By picking subsequences (without relabeling),
	we can assume that 
	the sequences
	$\seq{\hat\alpha_{k,i}}_k$,
	$\seq{\hat\beta_{k,j}}_k$
	are convergent with limits $\alpha_i \ge 0$, $\beta_j \in \R$.
	This shows
	\begin{equation*}
		\zeta_k
		=
		\xi_k
		-
		\sum_{i \in \AA(x)} \hat\alpha_{k,i} g_i - \sum_{j = 1}^m \hat\beta_{k,j} h_j
		\to
		\xi
		-
		\sum_{i \in \AA(x)} \alpha_{i} g_i - \sum_{j = 1}^m \beta_{j} h_j
		=:
		\zeta
	\end{equation*}
	with convergence in $L^{p'}(\mu)$.
	This implies $\zeta \in \NN_K(x)$ and, consequently,
	$\xi \in \NN_K(x) + \NN_P(x)$.

	The case $p = 1$ is slightly more difficult,
	since we have to work with the weak-$\star$ topology
	in $L^{p'}(\mu) = L^\infty(\mu)$.
	Thus, we have to consider weak-$\star$ convergent nets
	instead of sequences.
	By using the Krein--Šmulian theorem, \cite[Theorem~V.12.1]{Conway1985}, it is sufficient
	to consider a bounded net $\seq{\xi_k}_{k \in \mathcal{K}}$
	with
	\begin{equation*}
		\xi_k
		=
		\zeta_k + \sum_{i \in \AA(x)} \alpha_{k,i} g_i + \sum_{j = 1}^m \beta_{k,j} h_j,
	\end{equation*}
	where
	$\zeta_k \in \NN_K(x)$, $\alpha_{k,i} \ge 0$ and $\beta_{k,j} \in \R$,
	and
	$\seq{\xi_k}_{k \in \mathcal{K}} \weaklystar \xi$.
	Here, $\mathcal{K}$ is a directed set.
	We can argue as above to arrive at \eqref{eq:slater_conv834_bounds}.
	Now, we have to take an extra step,
	since convergent nets with values in $\R$
	do not need to be bounded.
	However, there exists $k_0 \in \mathcal{K}$ with
	\begin{equation*}
		0
		\ge
		\dual{\zeta_k}{\hat x - x}_{L^p(\mu)}
		\ge
		\dual{\xi}{\hat x - x}_{L^p(\mu)}
		-
		1
		\qquad\forall k \in \mathcal{K}, k \ge k_0
		.
	\end{equation*}
	Thus, we consider subnets
	indexed by
	\begin{equation*}
		\mathcal{K}' :=
		\set{
			k \in \mathcal{K}
			\given
			k \ge k_0
		}
		.
	\end{equation*}
	This is a directed set (with the order inherited from $\mathcal{K}$)
	and it is easy to check that
	$\seq{\xi_k}_{k \in \mathcal{K}'}$
	is a subnet of
	$\seq{\xi_k}_{k \in \mathcal{K}}$.
	Now, we can continue as above
	with some obvious changes.
\end{proof}
Together with \cref{lem:poly_cones},
we get the remarkable formula
\begin{equation}
	\label{eq:norm_K_P}
	\NN_{K \cap P}(x) = \NN_K(x) + \NN_P(x)
	,
\end{equation}
if we require the existence of a Slater point of \eqref{eq:linear}.
In particular,
by combining this result with \eqref{eq:opt_con_lin}
and
with the formulas for the normal cones,
we arrive at the following first-order condition.

\begin{theorem}
	\label{thm:foc_linear}
	Let $\bar x \in F$ be a locally optimal solution of \eqref{eq:linear}.
	Further, we assume the existence of a Slater point $\hat x$ of \eqref{eq:linear},
	cf.\ \cref{def:slater_point}.
	Then, there exist
	$\zeta \in L^{p'}(\mu)$,
	$\alpha_i \ge 0$, $i \in \AA(\bar x)$
	and
	$\beta_j \in \R$, $j = 1,\ldots, m$
	such that
	\begin{align*}
		f'(\bar x)
		+
		\zeta
		+
		\sum_{i \in \AA(\bar x)} \alpha_i g_i
		+
		\sum_{j = 1}^m \beta_j h_j
		&=
		0,
		\\
		\zeta &\le 0 \; \text{$\mu$-a.e.\ on $\set{x_a = \bar x < x_b}$},
		\\
		\zeta &= 0 \; \text{$\mu$-a.e.\ on $\set{x_a < \bar x < x_b}$},
		\\
		\zeta &\ge 0 \; \text{$\mu$-a.e.\ on $\set{x_a < \bar x = x_b}$}
		.
	\end{align*}
	Note that the multiplier $\zeta$
	can be split as
	$\zeta = \max(\zeta, 0) - \max(-\zeta, 0) =: \zeta_b - \zeta_a$
	and
	\begin{equation*}
		\zeta_a, \zeta_b \ge 0 \text{ on $\Omega$},
		\quad
		\zeta_a = 0 \text{ on $\set{x_a < \bar x}$}
		\quad
		\zeta_b = 0 \text{ on $\set{\bar x < x_b}$}
	\end{equation*}
	are satisfied $\mu$-a.e.\ by $\zeta_a, \zeta_b$.
\end{theorem}
Note that the Slater point $\hat x$ does not need to
satisfy the finitely many linear inequalities involving $g_i$
in a strict sense,
since the associated normal cone is automatically closed,
cf.\ \cite[Proposition~2.41]{BonnansShapiro2000}.

Finally, we give a small remark
which demonstrates that the Slater point
can ignore all ``real'' inequalities,
which cannot be reduced to equalities.
\begin{remark}
	\label{rem:real_inequalities}
	Let $I \subset \set{1,\ldots,n}$ be a subset of the inequality constraints
	which can be strictly satisfied, i.e.,
	there exists $\mathring x \in L^p(\mu)$
	with
	\begin{align*}
		x_a \le \mathring x &\le x_b \quad \text{$\mu$-a.e.},
		&
		\dual{g_i}{\mathring x}_{L^p(\mu)} &\le a_i \quad \forall i \in \set{1,\ldots,n} \setminus I,
		\\
		\dual{g_i}{\mathring x}_{L^p(\mu)} &< a_i \quad \forall i \in I,
		&
		\dual{h_j}{\mathring x}_{L^p(\mu)} &= b_j \quad \forall j = 1,\ldots,m.
	\end{align*}
	Further, we require the existence of $\tilde x \in L^p(\mu)$
	with
	\begin{align*}
		x_a < \tilde x &< x_b \quad \text{$\mu$-a.e.},
		&
		\dual{g_i}{\tilde x}_{L^p(\mu)} &\le a_i \quad \forall i \in \set{1,\ldots,n} \setminus I,
		\\
		&
		&
		\dual{h_j}{\tilde x}_{L^p(\mu)} &= b_j \quad \forall j = 1,\ldots,m.
	\end{align*}
	Then, it is easy to check that $(1-\varepsilon) \mathring x + \varepsilon \tilde x$
	is a Slater point for all $\varepsilon > 0$ small enough.
	Note that the point $\tilde x$ ignores all inequalities from $I$.
\end{remark}

\section{Necessity of Slater points}
\label{sec:slater_nec}
In this section, we are going to show
that the existence of a Slater point
is necessary for the existence of multipliers
in certain situations.
We continue to use the setting of \cref{sec:linear}.

First, we recall
the important result
that a finite linear system of equations and inequalities
can be reformulated equivalently
such that MFCQ holds.
\begin{lemma}
	\label{lem:mfcq_generic}
	There exists
	$\tilde x \in L^p(\mu)$,
	$\tilde n, \tilde m \in \N$,
	$\tilde g_i, \tilde h_j \in L^{p'}(\mu)$
	and
	$\tilde a_i, \tilde b_j \in \R$
	for
	$i = 1,\ldots,\tilde n$,
	$j = 1,\ldots,\tilde m$,
	such that
	\begin{enumerate}[label=(\roman*)]
		\item\label{lem:mfcq_generic:1}
			the linear systems
			\begin{align*}
				\dual{g_i}{x}_{L^p(\mu)} &\le a_i \quad \forall i = 1,\ldots, n,
				&
				\dual{h_j}{x}_{L^p(\mu)} &= b_j \quad \forall j = 1,\ldots, m
			\end{align*}
			and
			\begin{align*}
				\dual{\tilde g_i}{x}_{L^p(\mu)} &\le \tilde a_i \quad \forall i = 1,\ldots, \tilde n,
				&
				\dual{\tilde h_j}{x}_{L^p(\mu)} &= \tilde b_j \quad \forall j = 1,\ldots, \tilde m
			\end{align*}
			are equivalent for all $x \in L^p(\mu)$,
		\item\label{lem:mfcq_generic:2}
			the vectors $\tilde h_1, \ldots, \tilde h_{\tilde m}$
			are linearly independent,
		\item\label{lem:mfcq_generic:3}
			we have
			\begin{align*}
				\dual{\tilde g_i}{\tilde x}_{L^p(\mu)} &< \tilde a_i \quad \forall i = 1,\ldots, \tilde n,
				&
				\dual{\tilde h_j}{\tilde x}_{L^p(\mu)} &= \tilde b_j \quad \forall j = 1,\ldots, \tilde m.
			\end{align*}
	\end{enumerate}
\end{lemma}
\begin{proof}
	The conversion of our linear system can be performed in two steps.
	First, we convert all linear inequalities which cannot be strictly satisfied by feasible points
	into equalities.
	This ensures the existence of $\tilde x$ satisfying \ref{lem:mfcq_generic:3}.
	Second, we reduce the system of linear equalities to a maximal linear independent subset.
	This yields \ref{lem:mfcq_generic:2}.
	Finally, it is easy to see that both steps
	do not change the set of feasible points,
	thus \ref{lem:mfcq_generic:1} holds.
\end{proof}

\begin{lemma}
	\label{lem:density}
	Suppose that
	$x_a < x_b$ holds $\mu$-a.e.
	Then, a dense subset of $K$ is given by the set
	$\hat K := \set{x \in L^p(\mu) \given x_a < x < x_b \; \text{$\mu$-a.e.}}$.
\end{lemma}
\begin{proof}
	In the case that $K$ is empty, there is nothing to show.
	Otherwise, we choose some $\bar x \in K$.
	We first show that $\hat K$ is nonempty.
	Since $\mu$ is assumed to be $\sigma$-finite,
	it is easy to check that there exists a function $w \in L^p(\mu)$
	with $w > 0$.
	Then,
	a pointwise discussion shows that
	\begin{equation*}
		\hat x
		:=
		\bar x
		+
		\chi_{\set{\bar x = x_a}} \max\set{w, (x_b - \bar x)/2}
		-
		\chi_{\set{\bar x = x_b}} \max\set{w, (\bar x - x_a)/2}
	\end{equation*}
	lies strictly between $x_a$ and $x_b$
	and $\hat x \in L^p(\mu)$ follows from $\bar x, w \in L^p(\mu)$.
	Thus, $\hat x \in \hat K$.

	Finally, for an arbitrary $x \in K$,
	the sequence $\seq{ (1 - 2^{-n}) x + 2^{-n} \hat x }_{n \in \N}$
	belongs to $\hat K$ and converges towards $x$ in $L^p(\mu)$.
\end{proof}

\begin{lemma}
	\label{lem:no_slater}
	We suppose that $x_a < x_b$ and that the system \eqref{eq:linear} does not possess a Slater point.
	Then we have $\NN_P(\bar x) \cap -\NN_K(\bar x) \ne \set{0}$
	for all feasible points $\bar x$ of \eqref{eq:linear}.
\end{lemma}
\begin{proof}
	In the proof, we use the modified linear system from \cref{lem:mfcq_generic}.
	By
	\begin{equation*}
		\tilde\AA(\bar x)
		:=
		\set{ i \in \set{1,\ldots,\tilde n} \given \dual{\tilde g_i}{\bar x}_{L^p(\mu)} = \tilde a_i }
	\end{equation*}
	we denote the active set of $\bar x$
	and by $\tilde n(\bar x)$ we denote the cardinality of $\tilde\AA(\bar x)$.
	The assumption implies that
	\begin{align*}
		\dual{\tilde g_i}{x}_{L^p(\mu)} &\le \tilde a_i \quad \forall i \in \tilde\AA(\bar x),
		&
		\dual{\tilde h_j}{x}_{L^p(\mu)} &= \tilde b_j \quad \forall j = 1,\ldots, \tilde m
	\end{align*}
	is violated for all $x \in L^p(\mu)$ with $x_a < x < x_b$,
	since otherwise $(1-t) \bar x + t x$ would be a Slater point for $t \in (0,1)$ small enough.

	Thus, the sets
	\begin{align*}
		A &:=
		(\tilde a_i, \tilde b_j)
		_{\substack{i \in \tilde A(\bar x),\\j = 1,\ldots,\tilde m}}
		+ \R_-^{\tilde n(\bar x)} \times \set{0}^{\tilde m}
		,\\
		B &:=
		\set*{
			(\dual{\tilde g_i}{x}_{L^p(\mu)}, \dual{\tilde h_j}{x}_{L^p(\mu)})
			_{\substack{i \in \tilde A(\bar x),\\j = 1,\ldots,\tilde m}}
			% (\tilde g(x), \tilde h(x))
			\given
			x \in L^p(\mu),\; x_a < x < x_b
		}
	\end{align*}
	are disjoint, convex subsets of $\R^{\tilde n(\bar x) + \tilde m}$.
	A finite-dimensional separation theorem
	yields the existence of a vector $(\lambda, \mu) \in \R^{\tilde n(\bar x) + \tilde m} \setminus \set{0}$
	such that
	\begin{equation*}
		\lambda^\top v_1 + \mu^\top v_2
		\le
		\lambda^\top w_1 + \mu^\top w_2
		\qquad
		\forall
		(v_1, v_2) \in A,
		(w_1, w_2) \in B.
	\end{equation*}
	Since $v_1$ can be arbitrarily small,
	this gives $\lambda \ge 0$.
	Now, we observe that the vector
	\begin{equation*}
		(\dual{\tilde g_i}{\bar x}_{L^p(\mu)}, \dual{\tilde h_j}{\bar x}_{L^p(\mu)})
		_{\substack{i \in \tilde A(\bar x),\\j = 1,\ldots,\tilde m}}
	\end{equation*}
	belongs to $A$.
	Thus, the above inequality yields
	\begin{align*}
		0
		&\le
		\sum_{i \in \tilde\AA(\bar x)} \lambda_i \dual{\tilde g_i}{x - \bar x}_{L^p(\mu)}
		+
		\sum_{j = 1}^{\tilde m} \mu_j \dual{\tilde h_j}{x - \bar x}_{L^p(\mu)}
		\\&
		=
		\dual*{
			\sum_{i \in \tilde\AA(\bar x)} \lambda_i \tilde g_i
			+
			\sum_{j = 1}^{\tilde m} \mu_j \tilde h_j
		}{x - \bar x}_{L^p(\mu)}
	\end{align*}
	for all $x \in L^p(\mu)$ with $x_a < x < x_b$.
	Due to \cref{lem:density},
	the above inequality continues to hold for all $x \in K$
	and, thus,
	\begin{equation*}
		\zeta
		:=
		\sum_{i \in \tilde\AA(\bar x)} \lambda_i \tilde g_i
		+
		\sum_{j = 1}^{\tilde m} \mu_j \tilde h_j
		\in
		-\NN_K(\bar x).
	\end{equation*}
	Since the modified linear system from \cref{lem:mfcq_generic}
	still describes the polytope $P$,
	the sign condition on $\lambda$ and \cref{lem:finite_restrictions}
	show that $\zeta \in \NN_P(\bar x)$.
	It remains to check $\zeta \ne 0$.
	In case $\lambda \ne 0$,
	we use $\tilde x$ from \cref{lem:mfcq_generic}
	to obtain
	\begin{align*}
		\dual{\zeta}{\bar x - \tilde x}
		&=
		\sum_{i \in \tilde\AA(\bar x)} \lambda_i \dual{\tilde g_i}{\bar x - \tilde x}_{L^p(\mu)}
		+
		\sum_{j = 1}^{\tilde m} \mu_j \dual{\tilde h_j}{\bar x - \tilde x}_{L^p(\mu)}
		\\&
		=
		\sum_{i \in \tilde\AA(\bar x)} \lambda_i \dual{\tilde g_i}{\bar x - \tilde x}_{L^p(\mu)}
		>
		0
	\end{align*}
	since
	$\dual{\tilde g_i}{\bar x - \tilde x}_{L^p(\mu)} > 0$,
	$\lambda \ge 0$ and $\lambda \ne 0$.
	Hence, $\zeta$ cannot vanish.
	In the other case, we have $\lambda = 0$ and $\mu \ne 0$.
	By construction, the vectors $\tilde h_j$, $j = 1,\ldots,\tilde m$,
	are linearly independent and, again, we get that $\zeta \ne 0$.
\end{proof}

\begin{lemma}
	\label{lem:closure_normal}
	Let the assumptions of \cref{lem:no_slater}
	be satisfied.
	For a feasible point $\bar x$ of \eqref{eq:linear},
	let $\zeta \in \NN_P(\bar x) \cap - \NN_K(\bar x)$ with $\zeta \ne 0$
	be given.
	Then,
	\begin{equation*}
		\set*{
			\xi \in L^{p'}(\mu)
			\given
			\xi = 0 \text{ on $\set{\zeta = 0}$}
		}
		\subset
		\cl_\star\parens*{
			\NN_K(\bar x) + \NN_P(\bar x)
		}.
	\end{equation*}
\end{lemma}
\begin{proof}
	First, we note that $-\zeta \in \NN_K(\bar x)$
	yields
	$\bar x = x_b$ on $\set{\zeta < 0}$
	and
	$\bar x = x_a$ on $\set{\zeta > 0}$,
	see \cref{lem:normal_cone}.

	Let $\xi \in L^{p'}(\mu)$ with $\xi = 0$ on $\set{\zeta = 0}$ be given.
	We define the sequence $\seq{\xi_k}_{k \in \N} \subset L^{p'}(\mu)$
	via
	\begin{equation*}
		\xi_k(\omega) :=
		\begin{cases}
			\min\set{\xi(\omega), k \zeta(\omega) } & \text{if } \zeta(\omega) > 0, \\
			\max\set{\xi(\omega), k \zeta(\omega) } & \text{if } \zeta(\omega) < 0, \\
			0 & \text{else}.
		\end{cases}
	\end{equation*}
	Via a distinction by cases, we get
	\begin{align*}
		\xi_k - k \zeta &\le 0 \quad \text{on } \set{x_a = \bar x < x_b}, \\
		\xi_k - k \zeta & =  0 \quad \text{on } \set{x_a < \bar x < x_b}, \\
		\xi_k - k \zeta &\ge 0 \quad \text{on } \set{x_a < \bar x = x_b},
	\end{align*}
	i.e., $\xi_k - k\zeta \in \NN_K(\bar x)$, see again \cref{lem:normal_cone}.
	Thus,
	\begin{equation*}
		\xi_k = k \zeta + (\xi_k - k \zeta) \in \NN_P(\bar x) + \NN_K(\bar x).
	\end{equation*}
	Finally,
	the definition of $\xi_k$
	directly yields
	the pointwise $\mu$-a.e.\ convergence
	of $\xi_k$ towards $\xi$.
	Moreover, $\abs{\xi_k} \le \abs{\xi}$
	$\mu$-a.e.
	Thus,
	we can invoke Lebesgue's dominated convergence theorem
	to show that $\xi_k$ converges
	strongly (in case $p' < \infty$)
	or
	weak-$\star$ (in case $p' = \infty$)
	towards $\xi$.
	This yields the claim.
\end{proof}

Now, we are in position to prove the main theorem of this section.
\begin{theorem}
	\label{thm:no_slater_not_closed}
	We suppose that $x_a < x_b$ and that the system \eqref{eq:linear} does not possess a Slater point.
	Further, we assume that $p > 1$ (i.e., $p' < \infty$)
	and that the measure $\mu$ is nonatomic.
	Then, $\NN_K(\bar x) + \NN_P(\bar x)$
	is not closed in $L^{p'}(\mu)$
	for all feasible points $\bar x$ of \eqref{eq:linear}.
\end{theorem}
\begin{proof}
	Let $\zeta \in \NN_P(\bar x) \cap - \NN_K(\bar x)$ with $\zeta \ne 0$
	be given.
	We will construct a function $\xi \in L^{p'}(\mu)$
	with $\xi = 0$ on $\set{\zeta = 0}$
	and $\xi \not\in \NN_K(\bar x) + \NN_P(\bar x)$.
	Then, \cref{lem:closure_normal} yields the assertion.

	Since $\zeta$ does not vanish,
	one of the sets $\set{\zeta > 0}$ or $\set{\zeta < 0}$
	must have a positive measure.
	For simplicity, we assume $\mu(\set{\zeta > 0}) > 0$,
	in the other case we can argue similarly.
	Note that $\set{\zeta > 0} \subset \set{x_a = \bar x < x_b}$.

	We define an auxiliary function
	$\vartheta := \sum_{i = 1}^n \abs{g_i} + \sum_{j = 1}^m \abs{h_j}$
	and we note that every function from
	$\NN_K(\bar x) + \NN_P(\bar x)$
	is bounded from above by a scalar multiple of $\vartheta$
	on $\set{\zeta > 0}$,
	since the functions in $\NN_K(\bar x)$ are nonpositive on $\set{\zeta > 0}$.
	Since $\vartheta$ is real valued,
	we can find a subset $M \subset \set{\zeta > 0}$ of positive measure
	such that $\vartheta$ is bounded from above on $M$.
	Consequently, all functions in
	$\NN_K(\bar x) + \NN_P(\bar x)$
	are bounded from above on $M$.
	However, since $M$ has positive measure, since $p' < \infty$ and since $\mu$ is nonatomic,
	we can construct a nonnegative function $\xi \in L^{p'}(\mu) \setminus L^\infty(\mu)$,
	which vanishes outside of $M$.
	Now, \cref{lem:closure_normal} implies that
	$\xi \in \cl_\star\parens*{ \NN_K(\bar x) + \NN_P(\bar x) }$,
	while
	the boundedness of $\vartheta$ on $M$ and $\xi \not\in L^\infty(\mu)$
	gives
	$\xi \not\in \NN_K(\bar x) + \NN_P(\bar x)$.
	Thus, $\NN_K(\bar x) + \NN_P(\bar x)$ cannot be closed.
\end{proof}
The absence of Slater points also implies the nonexistence of multipliers.
\begin{corollary}
	\label{cor:no_slater_no_multiplier}
	We suppose that $x_a < x_b$ and that the system \eqref{eq:linear} does not possess a Slater point.
	Further, we assume that $p > 1$ (i.e., $p' < \infty$)
	and that the measure $\mu$ is nonatomic.
	Then, for each feasible point $\bar x$ of \eqref{eq:linear},
	there exists a linear functional $F \in L^{p'}(\mu)$
	such that $\bar x$ minimizes $F$ on $K \cap P$
	but no Lagrange multiplier exists.
\end{corollary}
\begin{proof}
	It is enough to choose a functional
	$F \in \NN_{K \cap P}(\bar x) = \cl_\star\parens*{ \NN_K(\bar x) + \NN_P(\bar x) }$
	satisfying $F \not\in \parens*{ \NN_K(\bar x) + \NN_P(\bar x) }$.
\end{proof}

The contrapositive of this corollary is also interesting.
It yields that, if we fix a feasible point $\bar x$
and if we obtain Lagrange multipliers
for all functionals which are minimized by $\bar x$ on $K \cap P$,
then there exists a Slater point.

\section{Nonlinear Constraints}
\label{sec:nonlinear}
In this section we consider the situation with additional nonlinear constraints.
For simplicity, we restrict ourselves to inequality constraints only.
That is, we study the problem
\begin{equation}
	\label{eq:nonlinear}
	\begin{aligned}
		\text{Minimize} \quad & f(x)\\
		\text{w.r.t.} \quad & x \in L^p(\mu) \\
		\text{s.t.} \quad
		&
		\begin{aligned}[t]
			x_a \le x &\le x_b && \text{$\mu$-a.e.\ on $\Omega$}, \\
			\dual{g_i}{x}_{L^p(\mu)} &\le a_i && \forall i = 1,\ldots, n,\\
			\dual{h_j}{x}_{L^p(\mu)} &= b_j && \forall j = 1,\ldots, m, \\
			G_i(x) &\le 0 && \forall i = 1,\ldots, N.
		\end{aligned}
	\end{aligned}
\end{equation}
Here, the data is as in \cref{sec:linear}
and, additionally, $N \in \N$
and
$G_i \colon L^p(\mu) \to \R$ are continuously Fréchet differentiable.
Similarly to \cref{sec:linear},
we require $G'(x) \in L^{p'}(\mu) = L^1(\mu)$
in the case $p = \infty$
for all feasible points $x$ of \eqref{eq:linear}.
In addition to the notation from \cref{sec:linear},
we define
\begin{equation*}
	Q := \set{x \in L^p(\mu) \given G_i(x) \le 0 \quad \forall i=1,\ldots,N}
\end{equation*}
and the feasible set of \eqref{eq:nonlinear}
is denoted by
\begin{equation*}
	E := F \cap Q = K \cap P \cap Q.
\end{equation*}
In order to linearize the constraints,
we utilize the constraint qualification
by
Robinson--Zowe--Kurcyusz (RZKCQ).
\begin{lemma}
	\label{lem:nonlin_RZK}
	Let $x \in E$ be a feasible point and we assume the existence of
	a point $\mathring x \in L^p(\mu)$ such that
	\begin{align*}
		x_a \le \mathring x &\le x_b \quad \text{$\mu$-a.e.\ on $\Omega$},
		&
		\dual{g_i}{\mathring x}_{L^p(\mu)} &\le a_i \quad \forall i = 1,\ldots, n,\\
		G_i'(x)(\mathring x - x) &< 0 \quad \forall i \in \BB(x),
		&
		\dual{h_j}{\mathring x}_{L^p(\mu)} &= b_j \quad \forall j = 1,\ldots, m,
	\end{align*}
	where
	\begin{equation*}
		\BB(x) := \set{ i \in \set{1,\ldots,N} \given G_i(x) = 0}
	\end{equation*}
	is the active set of the nonlinear constraints.
	Then,
	\begin{equation*}
		\TT_E(x)
		=
		\TT_{K \cap P}(x)
		\cap
		\TT_Q(x)
		=
		\TT_{K \cap P}(x)
		\cap
		\set{ h \in L^p(\mu) \given G'(x)h \le 0 \; \forall i \in \BB(x)}.
	\end{equation*}
\end{lemma}
\begin{proof}
	In order to apply the RZKCQ, we write the constraints in \eqref{eq:nonlinear}
	as
	\begin{equation*}
		x \in F, \qquad
		G(x) \in R,
	\end{equation*}
	where $G \colon L^p(\mu) \to \R^N$
	is the vector function with components $G_i$
	and $R = (-\infty,0]^N$.
	Since the interior of $R$ is nonempty,
	the RZKCQ can be formulated as
	\begin{equation*}
		\exists x' \in F :
		G(x) + G'(x) (x' - x) \in \operatorname{int}(R),
	\end{equation*}
	see \cite[(2.196)]{BonnansShapiro2000}.
	We can check that this condition is satisfied with $x' = (1-\varepsilon)x + \varepsilon \mathring x$
	for $\varepsilon > 0$ small enough.
	Thus, \cite[Corollary~2.91]{BonnansShapiro2000}
	yields the claim.
\end{proof}
We mention that RZKCQ also yields
a formula for the normal cone of $E$ if we define this cone
in the dual space of $L^p(\mu)$.
In case $p = \infty$, however, this is not compatible
with our \cref{def:polar,def:normal}.
Therefore, we provide the following lemma.

\begin{lemma}
	\label{lem:cones_Linfty}
	Let $A, B \subset L^\infty(\mu)$ be closed, convex cones
	such that $A - B = L^\infty(\mu)$.
	Further suppose that $B\polar \subset L^1(\mu)$ is contained
	in a finite-dimensional subspace of $L^1(\mu)$.
	Then, $A\polar + B\polar$ is closed in $L^1(\mu)$.
\end{lemma}
Note that the result \cite[Proposition~2.4.3]{Schirotzek2007}
(which does not need the assumption concerning finite dimensionality)
is not applicable here,
since, therein, the polars are defined to be subsets of $L^\infty(\mu)\dualspace$.
\begin{proof}
	From the generalized open mapping theorem \cite[Theorem~2.1]{ZoweKurcyusz1979},
	we get the existence of $C > 0$,
	such that for all $f \in L^\infty(\mu)$,
	there exist $g \in A$, $h \in B$ with $f = g - h$
	and $\max\set{\norm{g}_{L^\infty(\mu)}, \norm{h}_{L^\infty(\mu)}} \le C \norm{f}_{L^\infty(\mu)}$.

	Now let sequences
	$\seq{\zeta_k}_{k \in \N} \subset A\polar$
	and
	$\seq{\xi_k}_{k \in \N} \subset B\polar$
	be given
	such that $\zeta_k + \xi_k$ converges in $L^1(\mu)$.
	In particular, there exists $M \ge 0$ with $\norm{\zeta_k + \xi_k}_{L^1(\mu)} \le M$.
	Let $f \in L^\infty(\mu)$ with $\norm{f}_{L^\infty(\mu)} \le 1$
	be arbitrary.
	There exist $g \in A$, $h \in B$ with $f = g - h$
	and $\max\set{\norm{g}_{L^\infty(\mu)}, \norm{h}_{L^\infty(\mu)}} \le C$.
	Thus,
	\begin{align*}
		\dual{-\xi_k}{f}_{L^\infty(\mu)}
		&=
		\dual{\xi_k}{h - g}_{L^\infty(\mu)}
		\le
		\dual{\xi_k}{ - g}_{L^\infty(\mu)}
		\le
		\dual{\zeta_k + \xi_k}{ - g}_{L^\infty(\mu)}
		\\&
		\le
		\norm{\zeta_k + \xi_k}_{L^1(\mu)} \norm{-g}_{L^\infty(\mu)}
		\le
		C M
		.
	\end{align*}
	This bound holds for all $f$ with $\norm{f}_{L^\infty(\mu)} \le 1$,
	thus we get $\norm{\xi_k}_{L^1(\mu)} \le C M$.
	By assumption, $\seq{\xi_k}$ belongs to a finite-dimensional subspace,
	hence there exists a subsequence (without relabeling)
	such that $\xi_k \to \xi$ in $L^1(\mu)$.
	Consequently, $\zeta_k = (\zeta_k + \xi_k) - \xi_k$
	is also convergent.
	Since $A\polar$ and $B\polar$ are closed,
	this finishes the proof.
\end{proof}
By combining the previous lemmas,
we obtain a formula for the normal cone.
\begin{lemma}
	\label{lem:nonlin_normal_cone}
	Under the assumptions of \cref{lem:nonlin_RZK},
	we have
	\begin{equation*}
		\NN_E(x)
		=
		\NN_{K \cap P}(x)
		+
		\NN_Q(x),
	\end{equation*}
	where
	\begin{equation*}
		\NN_Q(x)
		=
		\set*{
			\sum_{i \in \BB(x)} \gamma_i G'_i(x)
			\given
			\gamma_i \ge 0, \; i \in \BB(x)
		}.
	\end{equation*}
\end{lemma}
\begin{proof}
	We are going to apply \cref{lem:cones_Linfty}
	with
	\begin{equation*}
		A = \TT_{K \cap P}(x)
		\quad\text{ and }\quad
		B = \TT_Q(x) = \set{ h \in L^p(\mu) \given G'(x)h \le 0 \; \forall i \in \BB(x)}
		.
	\end{equation*}
	For any $f \in L^\infty(\mu)$,
	we have
	\begin{equation*}
		G'(x) (-f + t (\mathring x - x)) < 0 \quad \forall i \in \BB(x)
	\end{equation*}
	if $t > 0$ is large enough,
	i.e., $-f + t(\mathring x - x) \in B$.
	Moreover, we have $t (\mathring x - x) \in \TT_{K \cap P}(x)$
	since $\mathring x \in K \cap P$ and $t > 0$.
	Thus,
	$f = t (\mathring x - x) - (-f + t (\mathring x - x)) \in A - B$.
	Thus, the claim follows from
	\begin{equation*}
		\set{ h \in L^p(\mu) \given G'(x)h \le 0 \; \forall i \in \BB(x)}\polar
		=
		\set*{
			\sum_{i \in \BB(x)} \gamma_i G'_i(x)
			\given
			\gamma_i \ge 0, \; i \in \BB(x)
		},
	\end{equation*}
	see \cite[Proposition~2.42]{BonnansShapiro2000},
	and from \cref{lem:cones_Linfty}.
\end{proof}
Finally,
it remains to utilize
\cref{lem:poly_cones,lem:slater}
to obtain the formula 
\eqref{eq:norm_K_P}
for $\NN_{K \cap P}(x)$.
Thus, we also need a Slater point as in \cref{def:slater_point}.
\begin{theorem}
	\label{thm:nonlin_RZK}
	Let $\bar x \in E$ be a feasible point and we assume the existence of
	$\tilde x \in L^p(\mu)$
	satisfying
	\begin{align*}
		x_a < \tilde x &< x_b \quad \text{$\mu$-a.e.\ on $\Omega$},
		&
		\dual{g_i}{\tilde x}_{L^p(\mu)} &\le a_i \quad \forall i = 1,\ldots, n,\\
		G_i'(\bar x)(\tilde x - \bar x) &< 0 \quad \forall i \in \BB(\bar x),
		&
		\dual{h_j}{\tilde x}_{L^p(\mu)} &= b_j \quad \forall j = 1,\ldots, m.
	\end{align*}
	where
	\begin{equation*}
		\BB(\bar x) := \set{ i \in \set{1,\ldots,N} \given G_i(\bar x) = 0}.
	\end{equation*}
	Then,
	\begin{equation*}
		\NN_E(\bar x)
		=
		\NN_K(\bar x)
		+
		\NN_P(\bar x)
		+
		\NN_Q(\bar x)
		.
	\end{equation*}
	Now, assume additionally that $x$ is a local minimizer of \eqref{eq:nonlinear}.
	Then, there exist
	$\zeta \in L^{p'}(\mu)$,
	$\alpha_i \ge 0$, $i \in \AA(\bar x)$,
	$\beta_j \in \R$, $j = 1,\ldots, m$
	and
	$\gamma_i \ge 0$, $i \in \BB(\bar x)$
	such that
	\begin{align*}
		f'(\bar x)
		+
		\zeta
		+
		\sum_{i \in \AA(\bar x)} \alpha_i g_i
		+
		\sum_{j = 1}^m \beta_j h_j
		+
		\sum_{i \in \BB(\bar x)} \gamma_i \gamma_i G'_i(\bar x)
		=
		0,&
		\\
		\zeta \le 0 \; \text{$\mu$-a.e.\ on $\set{x_a = \bar x < x_b}$}, &
		\\                                                               
		\zeta  =  0 \; \text{$\mu$-a.e.\ on $\set{x_a < \bar x < x_b}$}, &
		\\                                                               
		\zeta \ge 0 \; \text{$\mu$-a.e.\ on $\set{x_a < \bar x = x_b}$}. &
	\end{align*}
	The multiplier $\zeta$ can be split as in \cref{thm:foc_linear}.
\end{theorem}
\begin{proof}
	The result follows from \cref{lem:poly_cones,lem:slater,lem:nonlin_normal_cone}.
\end{proof}

On the first glance,
one might think that it would be
a weaker condition
to require the existence of
points $\hat x, \mathring x \in L^p(\mu)$
satisfying
\begin{align*}
	x_a < \hat x &< x_b \quad \text{$\mu$-a.e.\ on $\Omega$},
	&
	\dual{g_i}{\hat x}_{L^p(\mu)} &\le a_i \quad \forall i = 1,\ldots, n,\\
	&&
	\dual{h_j}{\hat x}_{L^p(\mu)} &= b_j \quad \forall j = 1,\ldots, m
	\shortintertext{and}
	x_a \le \mathring x &\le x_b \quad \text{$\mu$-a.e.\ on $\Omega$},
	&
	\dual{g_i}{\mathring x}_{L^p(\mu)} &\le a_i \quad \forall i = 1,\ldots, n,\\
	G_i'(\bar x)(\mathring x - \bar x) &< 0 \quad \forall i \in \BB(\bar x),
	&
	\dual{h_j}{\mathring x}_{L^p(\mu)} &= b_j \quad \forall j = 1,\ldots, m.
\end{align*}
Then, one can use $\mathring x$
in \cref{lem:nonlin_normal_cone}
to get $\NN_E(\bar x) = \NN_{K \cap P}(\bar x) + \NN_Q(\bar x)$
and afterwards $\hat x$
to get the representation of $\NN_{K \cap P}(\bar x)$.
However, similar to \cref{rem:real_inequalities},
we can set
$\tilde x = (1-\varepsilon)\mathring x + \varepsilon \hat x$.
If $\varepsilon > 0$ small enough,
$\tilde x$ satisfies the system in \cref{thm:nonlin_RZK}.

The conditions posed on the point $\tilde x$ from \cref{thm:nonlin_RZK}
are typically called
linearized Slater conditions.
As known from finite-dimensional optimization,
we need a strict inequality for nonlinear constraints
and a nonstrict inequality suffices for the (finite-dimensional)
linear constraints.
Moreover, the point $\tilde x$
is required to strictly satisfy
the linear (infinite-dimensional) box constraints.

\section{A counterexample}
\label{sec:counterex}

For illustration, we provide an example
which does not possess a Slater point.
We consider the problem \eqref{eq:linear}
with the following setting.
As a measure space, we use $\Omega = (0,1)$
equipped with the Lebesgue measure
and we take $p = p' = 2$.
We use the simple pointwise bounds
$x_a = 0$, $x_b = \infty$,
a single inequality constraint ($n = 1$)
$g_1 \equiv 1$, $a_1 = 0$,
and no equality constraints
($m = 0$).
In this situation, it is easy to check that
$\bar x = 0$ is the only feasible point
and the Slater condition w.r.t.\ $(K,P)$ is violated.
Here, $K$ and $P$ are as defined in \cref{sec:linear}.

Moreover, we have
\begin{align*}
	\NN_K(\bar x) &=
	\set{\zeta \in L^2(\mu) \given \zeta \le 0 \text{ $\mu$-a.e.\ on } (0,1)}
	&
	\NN_P(\bar x) &=
	\set{\alpha g_1 \given \alpha \ge 0}
	.
\end{align*}
For $\zeta \in L^2(\mu)$, it is easy to check
that $\zeta \in \NN_K(\bar x) + \NN_P(\bar x)$
if and only if the positive part $\max(\zeta,0)$
belongs to $L^\infty(\mu)$.
Thus, the closure of $\NN_K(\bar x) + \NN_P(\bar x)$
coincides with $L^2(\mu)$
and, in particular,
$\NN_K(\bar x) + \NN_P(\bar x)$
is not closed.
This is in accordance with
\cref{thm:no_slater_not_closed},
which says that
the assertion of
\cref{lem:slater}
fails in absence of a Slater point
and

Finally, we choose a function 
$z \in L^2(\mu)$ which is (essentially) unbounded from below
and we set $f(x) := \dual{z}{x}_{L^2(\mu)}$.
Then, $\bar x$ minimizes \eqref{eq:linear},
but the KKT conditions from \cref{thm:foc_linear}
cannot be satisfied.
The existence of such an objective $f$
is expected from \cref{cor:no_slater_no_multiplier}.

%%fakesection: Ack
\subsection*{Acknowledgement}
We would like to thank Fredi Tröltzsch
for bringing \cite{Troeltzsch1977}
to our attention.

%%fakesection: bibliography

\renewcommand*{\bibfont}{\small}
\printbibliography

\end{document}